%% file: paper.tex
\newif\ifarxiv
\arxivtrue
\pdfoutput=1
\RequirePackage[l2tabu, orthodox]{nag}

\documentclass[a4paper,final]{siamart220329}
\usepackage{enumitem}
\usepackage{amssymb}
\usepackage{tikz}
\usetikzlibrary{trees,perspective}
\usepackage{standalone}
\usepackage{xargs}
\usepackage{mathtools}
\usepackage{booktabs}
\usepackage{microtype}
\usepackage{cleveref}

\newlist{conjecturelist}{enumerate}{1}
\Crefname{conjecturelisti}{Conjecture}{Conjectures}
\crefname{conjecturelisti}{Conjecture}{Conjectures}

\let\div\relax
\DeclareMathOperator{\div}{div}
\DeclareMathOperator{\curl}{curl}
\DeclareMathOperator{\grad}{grad}
\newcommand{\beps}{{\boldsymbol \varepsilon}}
\newcommand{\bfz}{{\boldsymbol 0}}
\newcommand{\vt}{{\nu}}
\newcommand{\shift}{\sigma}
\newcommand{\kdg}{k}
\newcommand{\pmi}{n}
\newcommand{\sdiv}{{{\nabla\cdot} \,}}
\newcommand{\set}[2]{\left\lbrace #1 \; : \; #2 \right\rbrace}
\newcommand{\uu}{{u}}
\newcommand{\vv}{{v}}
\newcommand{\ww}{{w}}
\newcommand{\zz}{{z}}
\newcommand{\norm}[1]{\Vert #1 \Vert}
\newtheorem{conjecture}{Conjecture}
\newcommandx{\infsup}[4][1=u,2=v,3=V,4=V]{%
  \mathop{\mathrlap{\phantom{p}}\inf}_{\substack{{#1} \in {#3}\\ {#1} \neq 0}}%
  \;%
  \mathop{\sup}_{\substack{{#2} \in {#4} \\ {#2} \neq 0}}%
}

\headers{Two conjectures on the 3D Stokes complex}{P.~E.~Farrell, L.~Mitchell, AND L.~R.~Scott}

\title{Two conjectures on the Stokes complex in three dimensions on Freudenthal meshes}
\date{\today}
\author{
Patrick E.~Farrell\thanks{Mathematical Institute, University of Oxford, Oxford, England.
(\email{patrick.farrell@maths.ox.ac.uk})}
\and Lawrence Mitchell\thanks{NVIDIA Corporation, Santa Clara, CA, USA.\@
(\email{lmitchell@nvidia.com})}
\and L.~Ridgway Scott\thanks{Department of Computer Science, University of Chicago, Chicago, USA.\@
(\email{ridg@uchicago.edu}). PEF was funded by EPSRC grants EP/R029423/1 and EP/W026163/1.}
}

\begin{document}
\maketitle
\begin{AMS}
35Q30, 76D05, 65N30, 74S05, 76M10
\end{AMS}
\begin{abstract}
In recent years a great deal of attention has been paid to discretizations
of the incompressible Stokes equations that exactly
preserve the incompressibility constraint. These are of substantial
interest because these discretizations are pressure-robust, i.e.~the
error estimates for the velocity do not depend on the error in the pressure.
Similar considerations arise in nearly incompressible linear elastic solids.
Conforming discretizations with this property are now well understood
in two dimensions, but remain poorly understood in three dimensions.
In this work we state two conjectures on this subject. The first is
that the Scott--Vogelius element pair is inf-sup stable on uniform meshes
for velocity degree $k \ge 4$; the best result available in the literature is for $k \ge 6$.
The second is that there exists a stable space decomposition of the kernel of the divergence
for $k \ge 5$.
We present numerical evidence supporting our conjectures.
\end{abstract}

\section{Introduction}
We consider two closely related problems for a bounded Lipschitz domain $\Omega \subset \mathbb{R}^d$, $d \in \{2, 3\}$. The first is the incompressible Stokes
equations: given $f \in L^2(\Omega; \mathbb{R}^d)$, find the velocity $u: \Omega \to \mathbb{R}^d$ and pressure $p: \Omega \to \mathbb{R}$ such that
\begin{subequations}
\label{eq:stokes}
\begin{alignat}{3}
- \nabla \cdot \beps u + \nabla p &= f \quad && \text{ in } \Omega, \label{eq:stokes-momentum}\\
\nabla \cdot u &= 0 \quad && \text{ in } \Omega, \label{eq:stokes-incompressibility}\\
u &= 0 \quad && \text{ on } \partial \Omega,
\end{alignat}
\end{subequations}
where $\beps u = 1/2 \left(\nabla u + (\nabla u)^T\right)$ is the symmetric gradient of $u$. The second
is the Navier--Cauchy equation of linear elasticity: given $f \in L^2(\Omega; \mathbb{R}^d)$ and $\gamma > 0$, find the
displacement $u: \Omega \to \mathbb{R}^d$ that satisfies
\begin{subequations}
\label{eq:linear-elasticity}
\begin{alignat}{3}
- \nabla \cdot \beps u - \gamma \nabla \nabla \cdot u &= f \quad
  && \text{ in } \Omega, \label{eq:linear-elasticity-a}\\
u &= 0 \quad && \text{ on } \partial \Omega.
\end{alignat}
\end{subequations}
Here $\gamma = 2\lambda/\mu$, where $\lambda$ and $\mu$ are the Lam\'e parameters. As $\gamma \to \infty$,
the material is said to be nearly incompressible. The term $ - \nabla
\nabla \cdot u$ in \cref{eq:linear-elasticity-a} is connected
to the incompressibility constraint \cref{eq:stokes-incompressibility}; it arises 
in the Stokes momentum equation \cref{eq:stokes-momentum}
when employing an augmented Lagrangian approach~\cite{benzi2006}
to enforcing the divergence-zero constraint \cref{eq:stokes-incompressibility}.

When discretizing \cref{eq:stokes}, it is highly desirable to choose
spaces $V_h$ for the velocity and $\Pi_h$ for the pressure such that
all discretely divergence free functions are pointwise divergence free,
i.e.~the
incompressibility constraint \cref{eq:stokes-incompressibility} is
satisfied exactly on the discrete level~\cite{john2017}.
Achieving this is difficult; no element pair for exact enforcement is known that
is simultaneously inf-sup stable, low-order, conforming, has polynomial basis functions, and is effective on
general meshes. On simplicial grids with special mesh structure, it is possible to
use the conforming Scott--Vogelius finite element pair
$[\mathrm{CG}_k]^d$-$\mathrm{DG}_{k-1}$~\cite{lrsBIBbk,ScottVogeliusA} for $k \ge d$ (for Alfeld
splits~\cite{ref:QinThesis,ref:GuzmaNeilanBarycntr,fu2020exact,ref:zhang3DalfeldSplit})
or $k \ge d-1$ (for Powell--Sabin
splits~\cite{ref:GuzmanLisNeilanPowelSabin,ref:linearPowellSabinZhang,ref:quadraticPowellSabinTets,guzman2020exact}). The approach of
Guzm\'an and Neilan~\cite{ref:GuzmaNeilanBarycntr} is conforming, works for
arbitrary degree and on general meshes, but requires the use of
piecewise polynomial basis functions on each cell (instead of standard polynomials). Non-conforming divergence-free
discretizations are reviewed in John et al.~\cite[\S
4.4]{john2017}. Finally, another approach is to consider
the use of high-order discretizations, which are attractive for
their advantageous computational properties on modern architectures~\cite{ScottVogeliusA,lrsBIBbk,neilan2020stokes,zhang2011divergence}. Another alternative is to modify the right-hand side of the problem with an operator that maps discretely divergence-free test functions to exactly divergence-free ones~\cite{linke2012}.

In this paper we state two conjectures regarding the discretization and
multigrid solution of \cref{eq:stokes} and \cref{eq:linear-elasticity} on structured uniform tetrahedral meshes
(Freudenthal meshes~\cite{ref:Freudenthalulation}). We focus on Freudenthal meshes since some theoretical results
are known in this case. For concreteness we briefly describe the
Freudenthal triangulation of the unit cube in \cref{sec:freudenthal}.
The first is that the Scott--Vogelius element pair is inf-sup stable on
Freudenthal meshes
for $k \ge 4$; the best available result is that of Zhang, who
proved that the pair is stable for $k \ge 6$~\cite{zhang2011divergence}. We conjecture this on the
basis of numerical calculations of the inf-sup constant for varying $k$.
These rely on a new algorithm that can compute the inf-sup constant for elements that are divergence-free when an exact characterization of the pressure space is not known.
The second is that on the same meshes the subspace $Z_h \subset V_h$ of divergence-free
functions admits a local basis defined on the vertex-centred patches for
$k \ge 5$, and that the associated space decomposition is stable. This is significant
because identifying a local basis for the kernel of the divergence operator is crucial
for multigrid algorithms applied to \cref{eq:linear-elasticity}.
We conjecture this on the basis of observed $\gamma$-robustness of a
multigrid solution algorithm for \cref{eq:linear-elasticity}, for which
the local kernel decomposition is essential~\cite{schoberl1999}.
The existence of a local basis is known in three dimensions for Alfeld splits~\cite{fu2020exact} and Worsey--Farin splits~\cite{guzman2020exact}, but remains an open question for general meshes, and in particular for the Freudenthal meshes considered here.

\section{Inf-sup stability of the Scott--Vogelius element}

\begin{figure}
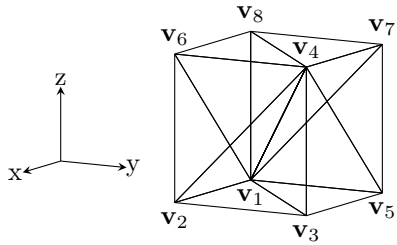

  \centering
  \includestandalone{freudenthal-subdivision}
  \caption{Unit cube for a structured mesh in three dimensions used for
    computational experiments. Given $N \in \mathbb{N}$, a $N \times N \times N$ mesh
    of cubes is generated, each of which is subdivided into 6 tetrahedra as shown. Vertex labels correspond to the vertices enumerated in \cref{sec:freudenthal}.}
  \label{fig:freudenthal-subdivision}
\end{figure}

The mixed formulation of \cref{eq:stokes} is to find 
$(u, p) \in V \times \Pi \coloneqq H^1_0(\Omega; \mathbb{R}^d) \times L^2_0(\Omega)$
such that
\begin{equation} \label{eq:stokes-weak}
(\beps u, \beps v)_{L^2} - (p, \sdiv v)_{L^2} - (q, \sdiv u)_{L^2} =
( f, v )_{L^2(\Omega)} \quad \text{ for all } (v, q) \in V \times \Pi.
\end{equation}
Define $( T , U )_{L^2(\Omega)} = (T, U)_{L^2}$ for any scalar-, 
vector-, or tensor-valued functions $T,U$.

The inf-sup condition determines whether or not a 
pair of spaces $V_h \subset V$, $\Pi_h \subset \Pi$ in a mixed finite element method provide a compatible discretization~\cite{lrsBIBgd}, since the symmetric gradient term is coercive on $V$.
In the context of \cref{eq:stokes-weak}, this condition encapsulates a constraint between the
divergence of the velocity space $V_h$ and the pressure space $\Pi_h$: there exists
$\beta > 0$ such that
\begin{equation} \label{eqn:infsupcon}
 \beta\, \sqrt{( q , q )_{L^2(\Omega)}} \leq
 \sup_{\substack{v\in V_h \\ v \neq 0}}\frac{( \sdiv v,q)_{L^2(\Omega)}}
{ \sqrt{( \nabla v , \nabla v )_{L^2(\Omega)}} }
\quad\forall q\in\Pi_h.
\end{equation}

In this section we make the following conjecture:
\begin{conjecture} \label{con:conjectureone}
Let $V_h$ be constructed with continuous Lagrange elements of degree $k$, and choose $\Pi_h = \sdiv V_h$. For $\kdg \ge 4$, the inf-sup condition \cref{eqn:infsupcon} holds on structured tetrahedral 
meshes of the form in \cref{fig:freudenthal-subdivision}, known as the
Freudenthal triangulation~\cite{ref:Freudenthalulation}, with a constant
that depends only on $\kdg$.
\end{conjecture}

We conjecture this on the basis of numerical computations with a new algorithm
for calculating the inf-sup constant, which we now describe.

\subsection{Computing the inf-sup constant}

There are several methods to estimate computationally
the inf-sup constant $\beta$ for various spaces and variational problems~\cite{ref:QinThesis,ref:RognesAutoInfSupConst,ref:ArnoldRognesInfSupStabil}.
In~\cite{ref:RayleighRitzInfSupConst}, approximation of the corresponding
Ladyzhenskaya inf-sup constant for the continuous problem is studied.
All of these are based on solving eigenproblems.

In~\cite{ref:RognesAutoInfSupConst}, an automated system ASCoT was introduced
for computing inf-sup constants.
It takes as input $V_h$ and $\Pi_h$.
However, when $\Pi_h=\sdiv V_h$, it is advantageous to exploit this structure
in the calculation of the inf-sup constant.
More importantly, the space $\sdiv V_h$ may not be known {\em a priori}, and
so a different algorithm is needed that does not require $\Pi_h$ as input.

A technique was developed by one of the authors in~\cite{lrsBIBih} that deals specifically with the 
$\Pi_h=\sdiv V_h$ case and casts this as an eigenvalue problem, as follows.
Define the bilinear forms
\begin{equation}\label{eqn:gradform}
a(u,v)= \int_\Omega \nabla u : \nabla v \,dx ,\qquad
b(v,q)= \int_\Omega (\sdiv v)q \,dx .
\end{equation}
First define
\begin{equation}
Z_h = \set{v \in V_h}{b(v, q) = 0 \;\forall q \in \Pi_h},
\end{equation}
which is the set of divergence-free functions in $V_h$ since $\Pi_h=\sdiv V_h$.
Define $\kappa$ by
\begin{equation}\label{eqn:bisbetaabndt}
\kappa=\min_{0\neq v\in V_h, \; v\perp_a Z_h} \frac{(\sdiv v,\sdiv v)_{L^2}}{a(v,v)}
      =\min_{0\neq v\in Z_h^\perp} \frac{(\sdiv v,\sdiv v)_{L^2}}{a(v,v)},
\end{equation}
where $v\perp_a Z_h$ means that $a(v,w)=0$ for all $w\in Z_h$
and 
\begin{equation}
Z_h^\perp=\set{v\in V_h}{a(v,w)=0 \;\forall w\in Z_h}.
\end{equation}
We recall the following lemma from~\cite[Lemma 26.1]{lrsBIBih}.

\begin{lemma}
Suppose that $\norm{v}_V=\sqrt{a(v,v)}$ and $\Pi_h=\sdiv V_h$. Then
\begin{equation}\label{eqn:fomdugook}
\beta = \infsup[q][v][\Pi_h][V_h] \frac{b(v,q)}{\norm{v}_V\norm{q}_{L^2}}
      = \infsup[q][v][\Pi_h][Z_h^\perp] \frac{b(v,q)}{\norm{v}_V\norm{q}_{L^2}}
      \geq\sqrt{\kappa}\geq \tfrac{1}{2} \beta,
\end{equation}
where $\kappa$ is defined in \cref{eqn:bisbetaabndt}.
\end{lemma}

This result holds also for the bilinear form 
\begin{equation}\label{eqn:epsformak}
 a(u,v)= ( \beps u , \beps v )_{L^2(\Omega)} ,
\end{equation}
with possibly different constants.

Computing $\kappa$ is equivalent to finding the smallest eigenvalue $\lambda$
of the following eigenproblem: find $0\neq u_h\in Z_h^\perp$ such that
\begin{equation}\label{eqn:farfrmsv}
(\sdiv u_h,\sdiv v_h)_{L^2(\Omega)}=\lambda\, a(u_h,v_h) \quad\forall v_h\in Z_h^\perp,
\end{equation}
which is equivalent to the Rayleigh quotient minimization \cref{eqn:bisbetaabndt}.
Note that $\kappa>0$ since $\kappa=0$ leads to the contradiction
$\sdiv u_h=0$, that is, $u_h\in Z_h\cap Z_h^\perp$.
Thus there are no spurious modes when $\Pi_h=\sdiv V_h$.

We can write \cref{eqn:farfrmsv} in operator form as
\begin{equation}\label{eqn:operfmsv}
B u_h=\lambda A u_h.
\end{equation}
The equation \cref{eqn:operfmsv} is a symmetric generalized eigenvalue problem~\cite{stewart2001matrix}, and its eigenvalues are all real.
There are many algorithms for solving symmetric generalized eigenvalue problems~\cite{golub2002inverse,sameh2000trace,stewart2001matrix}.
However, we do not have an explicit basis for the space $Z_h^\perp$, so we
choose to use matrix-free methods.
Here we focus on a simple method related to the power method.

Since $A$ is invertible on all of $V_h$, it would be attractive to utilize
an iteration in which we invert $A$, and not $B$.
We introduce a shift $\shift$:
\begin{equation}\label{eqn:shtrfmsv}
(B-\shift A)\uu_h=\lambda^\shift A \uu_h.
\end{equation}
If $\lambda_1\leq\lambda_2\leq\cdots\leq\lambda_N$ are the eigenvalues of 
the symmetric problem \cref{eqn:operfmsv}, that is, for the operator $A^{-1}B$,
and $\lambda^\shift_1\leq\lambda^\shift_2\leq\cdots\leq\lambda^\shift_N$
are the eigenvalues of the shifted problem \cref{eqn:shtrfmsv}, then 
$\lambda^\shift_i=\lambda_i-\shift$ for all $i=1,\dots,N$.
Moreover, $\lambda^\shift_i$ are the eigenvalues for the operator $A^{-1}B-\shift I$,
and the eigenvectors $\uu_h^i$ for $\lambda^\shift_i$ and $\lambda_i$ are the same,
for all $i=1,\dots,N$.

\subsection{Using the power method}

We can solve for certain eigenvalues and eigenvectors $\uu_h$ via 
the power method~\cite{lrsBIBgh} for \cref{eqn:shtrfmsv}, namely to find 
$\uu^\pmi\in Z_h^\perp$ such that
\begin{equation}\label{eqn:rayquosv}
\begin{split}
a(\uu^{\pmi},\vv_h) &= (\sdiv\hat\uu^{\pmi-1},\sdiv\vv_h)_{L^2(\Omega)} -\shift\, a(\hat\uu^{\pmi-1},\vv_h) 
\quad\forall\;\vv_h\in Z_h^\perp\\
\lambda^{\pmi+1}&=\frac{(\sdiv\uu^\pmi,\sdiv\uu^\pmi)_{L^2(\Omega)}}{a(\uu^\pmi,\uu^\pmi) }, \qquad
\hat\uu^{\pmi}=(\norm{\uu^{\pmi}}_V)^{-1} \uu^\pmi.
\end{split}
\end{equation}
The choices of initial iterate and termination criterion are discussed subsequently.

We now consider how to compute this, despite the fact that the space $Z_h^\perp$
is not explicitly known.
Suppose that $\hat\uu^0\in Z_h^\perp$ is given and that we solve for 
$\uu^{1}\in V_h$ via
$$
a(\uu^{1},\vv_h) = (\sdiv\hat\uu^{0},\sdiv\vv_h)_{L^2(\Omega)}
-\shift\, a(\hat\uu^{0},\vv_h) 
\quad\forall\; \vv_h\in V_h,
$$
which we can do since $a(\cdot,\cdot)$ is coercive on $V_h$.
Then for all $\vv_h\in Z_h$
$$
a(\uu^{1},\vv_h) = (\sdiv\hat\uu^{0},\sdiv\vv_h)_{L^2(\Omega)} -\shift\, a(\hat\uu^{0},\vv_h) = 0
$$
since both terms vanish.
Thus $\uu^{1}\in Z_h^\perp$.
Moving on, we can solve for $\uu^{\pmi}\in V_h$ via
\begin{equation}\label{eqn:howesolv}
a(\uu^{\pmi},\vv_h) = (\sdiv\hat\uu^{\pmi-1},\sdiv\vv_h)_{L^2(\Omega)}-\shift\, a(\hat\uu^{\pmi-1},\vv_h) 
\quad\forall\;\vv_h\in V_h.
\end{equation}
By induction, $a(\uu^\pmi,\vv_h)=0$ for all $\vv_h\in Z_h$, so that $\uu^{\pmi}\in Z_h^\perp$
for all $\pmi> 0$.

To start the process, we can solve for $\uu^{0}\in V_h$ via
\begin{equation}\label{eqn:starproc}
a(\uu^{0},\vv) = (\sdiv\ww,\sdiv\vv)_{L^2(\Omega)} \quad\forall\;\vv\in V_h,
\end{equation}
where we must pick some $\ww$ where $\sdiv\ww\neq 0$.
In practice, we chose
\begin{equation}\label{eqn:kaywhyin}
\ww=(\sin(\omega x),\cos(\omega y)),
\end{equation}
for $\omega \in \mathbb{Z}$,
although other initializations worked as well.
We then set 
\begin{equation}\label{eqn:sstedtroc}
\hat\uu^{0}=(\norm{\uu^{0}}_V)^{-1} \uu^0.
\end{equation}
Note that $a(\uu^{0},\vv) = 0$ for all $\vv\in Z_h$, which means that
$\uu^{0}\in Z_h^\perp$.
The division in \cref{eqn:sstedtroc} provides a natural check that $\uu^{0}\neq\bfz$.

\subsection{Eigenvalue bounds}

The eigenvalues $\lambda_i$ are bounded above since
$$
\norm{\sdiv\vv_h}_{L^2}^2 \leq C a(\vv_h,\vv_h) \quad\forall\;\vv_h\in V.
$$
This is obvious for the gradient form \cref{eqn:gradform}
with $C\leq d$, where $d$ is the dimension of $\Omega$.
For the $\beps$ form \cref{eqn:epsformak}, it follows by the Korn inequality~\cite[\S 11.2]{lrsBIBgd}.
A good estimate of $C$ can be obtained computationally since it is the largest
eigenvalue of \cref{eqn:farfrmsv}.
In all computational examples here, using the gradient form \cref{eqn:gradform},
 $C=1$, apparently due to the homogeneous Dirichlet boundary conditions in $V$.
Using the $\beps$ form \cref{eqn:epsformak}, which is equivalent to \cref{eqn:gradform} by
Korn's inequality, would change the computational results here by at most a constant factor.
When homogeneous Dirichlet boundary conditions are enforced, it is 
known~\cite[(13.12)]{lrsBIBih} that the forms \cref{eqn:epsformak} and \cref{eqn:gradform} 
produce identical values for divergence-free functions.

Thus we can shift by a constant $\shift$ independent of $h$ to ensure that
$\lambda^\shift_i<0$ for all $i$.
The algorithm \cref{eqn:rayquosv} is the power method for $A^{-1}B-\shift I$,
and it will generically converge to the eigenvector associated with
the most negative eigenvalue, provided $\shift>\tfrac{1}{2}\lambda_N$, and
so it will generically converge to $\lambda_1^\shift$.
Taking $\shift=0$ will converge to the largest eigenvalue $\lambda_N$, since $0 < \lambda_1$.
We found that taking $\shift=0.6$ gave acceptable convergence.

\subsection{Effect of round-off error}

It is essential to start the iteration \cref{eqn:rayquosv} with 
$\hat\uu^{0}\in Z_h^\perp$.
To see what goes wrong otherwise, write the algorithm in \cref{eqn:rayquosv} as
$$
\uu^{k+1}=c_k (A^{-1}B-\shift I)\uu^k,
$$
where $c_k=1/\norm{\uu^k}_V$.
Suppose that $\uu^0=\vv^0+\ww^0$ where $\vv^0\in Z_h^\perp$ and
$\ww^0\in Z_h$.
Then $\uu^k=\vv^k+\ww^k$ where $\vv^k\in Z_h^\perp$ and $\ww^k\in Z_h$, and
$$
\vv^{k+1}=c_k (A^{-1}B-\shift I)\vv^k, \qquad
\ww^{k+1}=c_k (-\shift I)\ww^k.
$$
The reason is that $A^{-1}B$ maps $Z_h^\perp$ into itself and $Z_h$ to
the zero vector.
Define $C_k=\prod_{i=0}^{k-1} c_i$.
Then
$$
\vv^{k}=C_k (A^{-1}B-\shift I)^k \vv^0, \qquad
\ww^{k}=C_k (-\shift I)^k \ww^0.
$$
In seeking the smallest eigenvalue $\lambda_1$, we will need to take $\shift$
as large as at least half of the largest eigenvalue $\lambda_N$.
Thus $\ww^k$ can become significant even if $\ww^0$ is on the order of
round-off error.
Once $\ww^k$ becomes dominant, the Rayleigh quotient in \cref{eqn:rayquosv}
defining $\lambda^k$ will go to zero.
Thus it is necessary to monitor the projection $\zz^k$ of $\uu^k$ onto $Z_h$ and
to project $\uu^k$ onto $Z_h^\perp$ (by simply subtracting $\zz^k$)
when the divergence-zero component is too large.
This may need to be done for $k=0$ as well.

The projection $\zz^k\in Z_h$ satisfies
$$
a(\zz^k,\vv)= a(\uu^k,\vv) \;\forall\;\vv\in Z_h.
$$
The projection $\zz^k$ can be computed via the iterated penalty method.
Some care is required in using the iterated penalty method to do this, since
it will be slow to converge if the inf-sup constant is small.
But this appears to work well in practice, albeit with a large number of 
iterations required when the inf-sup constant is smaller.

\subsection{Full algorithm}

We now summarize the full algorithm.
First we compute $\hat\uu^{0}$ via \cref{eqn:starproc} and \cref{eqn:sstedtroc}, 
where $\ww$ is given in \cref{eqn:kaywhyin}.
Next, we solve for $\uu^{\pmi}$ via \cref{eqn:howesolv}, $n\geq 1$.
Then we project $\uu^{\pmi}$ onto $Z_h$ using the iterated penalty
method~\cite[\S 13.1]{lrsBIBgd}.
To compute the projection of $\uu^{\pmi}$ onto $Z_h$, we solve for $\zz^{\ell},\;\ww^{\ell} \in V_h$ such that
\begin{equation}\label{eqn:fulrayquosv}
\begin{split}
&a(\zz_{\ell},\vv) + \rho (\sdiv\zz_{\ell},\sdiv\vv) = a(\uu^{\pmi},\vv)-(\sdiv\ww_{\ell},\sdiv\vv)
\quad\forall\;\vv\in V_h\\
&\ww_{\ell+1}=\ww_{\ell+1}+\rho \zz_{\ell},
\end{split}
\end{equation}
where we start with $\ww_{0}=0$. The parameter $\rho$ is the penalty parameter enforcing the incompressibility constraint; in our computations we set $\rho = 10^4$.
We terminate the iteration on $\ell$ when 
\begin{equation}\label{eqn:termdyquosv}
\norm{\sdiv\zz_{\ell}}_{L^2(\Omega)}\leq \tau,
\end{equation}
where we picked $\tau=10^{-14}$ in our computations.
Then $\zz_{\ell}\approx \Pi_{Z_h}\uu^{\pmi}$, and if 
$$
\norm{\nabla\zz_{\ell}}_{L^2(\Omega)}\geq \zeta \norm{\nabla\uu^{\pmi}}_{L^2(\Omega)},
$$
we update
$
\uu^{\pmi}=\uu^{\pmi}-\zz_{\ell}.
$
In our computations, we picked $\zeta=10^{-12}$.
Once $\uu^{\pmi}$ has been computed, we define
\begin{equation}\label{eqn:lasrayquosv}
\begin{split}
&\lambda^{\pmi+1}=\frac{(\sdiv\uu^\pmi,\sdiv\uu^\pmi)}{a(\uu^\pmi,\uu^\pmi) }, \qquad
\hat\uu^{\pmi}=(\norm{\uu^{\pmi}}_V)^{-1} \uu^\pmi.
\end{split}
\end{equation}
This iteration is continued while $|\lambda^{\pmi+1}-\lambda^{\pmi}|>\epsilon$, where 
$\epsilon$ is a pre-specified tolerance, with $\epsilon=10^{-8}$ in our computations.

\subsection{2D tests}

\begin{table}
\centering
\caption{Mesh restrictions for $V_h$ constructed with continuous piecewise polynomials of degree $\kdg$, and $\Pi_h = \sdiv V_h$. 
Here $d$ is the dimension of $\Omega$.
In the lowest-order case $k=1$ in two dimensions, the velocity approximation
is optimal order despite the fact that the inf-sup constant deteriorates.\label{tabl:divfreeOK}}
\begin{tabular}[t]{c|c|c|l}
\toprule
$d$ & $\kdg$ & inf-sup & mesh restrictions \\
\midrule
  2 & $\phantom{\geq}$ 1 & no, but & optimal velocity approximation on Malkus splits~\cite{malkus1984linear,ref:QinZhangCrossedTriangles}\\
  2 & $\phantom{\geq}$ 1 & yes & Powell--Sabin splits~\cite{ref:loworderWorseyFarin,ref:GuzmanLisNeilanPowelSabin,ref:linearPowellSabinZhang}\\
  2 & $\phantom{\geq}$ 2 & yes & some crossed triangles required~\cite{ref:QinThesis} 
or Alfeld splits~\cite{ref:GuzmaNeilanBarycntr,fu2020exact} \\
  2 & $\phantom{\geq}$ 3 & yes & 
under certain conditions~\cite{lrsBIBic}
\\
  2 & $\geq 4$ & yes & $p$-robust when no nearly singular vertices~\cite{lrsBIBbk,ScottVogeliusA,lrsBIBib,ainsworth2022} \\
\midrule
  3 & $\geq 1$ & yes & Worsey--Farin splits~\cite{ref:loworderWorseyFarin} \\
  3 & $\geq 3$ & yes & Alfeld splits~\cite{ref:zhang3DalfeldSplit,ref:GuzmaNeilanBarycntr,fu2020exact} \\
  3 & $\geq 6$ & yes & only one family of meshes known~\cite{zhang2011divergence}\\
\bottomrule
\end{tabular}
\end{table}

\begin{table}
\centering
\caption{Computation of inf-sup constants
on Malkus splits (Type II meshes) in two dimensions.
The shift parameter in \cref{eqn:shtrfmsv} was set to $\shift=0.6$. 
The mesh size $N$ refers to the $N \times N$ mesh of quadrilaterals before subdivision.
The parameter $\omega$ determines the initial eigenvector approximation defined 
in \cref{eqn:kaywhyin}.
Iterations were continued until the change in eigenvalue was less than $10^{-7}$.\label{tabl:lintwode}}
\begin{tabular}[t]{c|ccccc}
\toprule
degree $\kdg$ & $N$ & $\omega$ & iterations & inf-sup $\lambda_1$    & restarts \\
\midrule
1             & 5   & 5        & 47         & $4.08 \times 10^{-1}$  & 0        \\ 
1             & 5   & 10       & 40         & $4.08 \times 10^{-1}$  & 0        \\ 
1             & 5   & 100      & 51         & $4.08 \times 10^{-1}$  & 1        \\ 
1             & 10  & 10       & 258        & $1.13 \times 10^{-2}$  & 1        \\
1             & 20  & 10       & 642        & $2.98 \times 10^{-3}$  & 1        \\ 
\midrule
2             & 10  & 10       & 193        & $1.49 \times 10^{-1}$  & 13       \\
2             & 20  & 10       & 189        & $1.48 \times 10^{-1}$  & 14       \\
2             & 40  & 10       & 187        & $1.48 \times 10^{-1}$  & 19       \\
\bottomrule
\end{tabular}
\end{table}

\begin{figure}
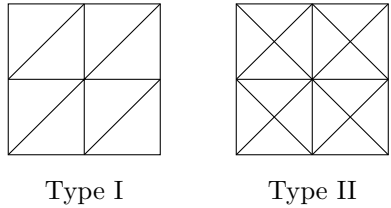

  \centering
  \includestandalone{twotypes}
  \caption{Two types of regular meshes in two dimensions.
    The Type II mesh is also known as the Malkus split.}
  \label{fig:twotypes}
\end{figure}

\begin{table}
\centering
\caption{Computation of inf-sup constants
on Type I meshes (see \cref{fig:twotypes}) in two dimensions. 
The shift parameter in \cref{eqn:shtrfmsv} was set to $\shift=0.6$.
The mesh size $N$ refers to the $N \times N$ mesh of quadrilaterals before subdivision.
The parameter $\omega$ determines the initial eigenvector approximation defined 
in \cref{eqn:kaywhyin}.
Iterations were continued until the change in eigenvalue was less than $10^{-7}$
for $k=3$ and $k=4$, but we had to continue iterations 
until the change in eigenvalue was less than $10^{-10}$ for $k=2$ to obtain a reliable second digit.\label{tabl:cubqwode}}
\begin{tabular}[t]{c|ccccc}
\toprule
degree $\kdg$ & $N$ & $\omega$ & iterations & inf-sup $\lambda_1$      & restarts \\
\midrule
2             & 8   & 10       & 13577      & $1.60 \times 10^{-3}$    & 29       \\ 
2             & 16  & 10       & 50690      & $4.08 \times 10^{-4}$    & 129      \\ 
2             & 32  & 10       & 19851      & $1.07 \times 10^{-4}$    & 19397    \\ 
2             & 64  & 10       & 54650      & $2.75 \times 10^{-5}$    & 54649    \\ 
\midrule
3             & 3   & 10       & 790        & $8.46 \times 10^{-3}$    & 11       \\ 
3             & 5   & 10       & 768        & $3.52 \times 10^{-3}$    & 8        \\ 
3             & 10  & 10       & 1858       & $9.50 \times 10^{-4}$    & 15       \\ 
\midrule
4             & 5   & 10       & 585        & $2.59 \times 10^{-2}$    & 12       \\ 
4             & 10  & 10       & 789        & $2.60 \times 10^{-2}$    & 17       \\ 
4             & 20  & 10       & 1265       & $2.60 \times 10^{-2}$    & 28       \\ 
\bottomrule
\end{tabular}
\end{table}

To verify the algorithm, we compare with known results.
We summarize some known results in two and three dimensions in \cref{tabl:divfreeOK}.
In \cref{tabl:lintwode}, we give results for Lagrange elements of degree $\kdg=1,2$
in two dimensions on Malkus splits (Type II meshes in the nomenclature of~\cite{lai2007spline}
as shown in \cref{fig:twotypes}), for various mesh sizes.
A Malkus split mesh starts with a quadrilateral subdivision and creates a triangulation
by subdividing each quadrilateral by adding the two diagonals connecting opposite vertices.
On Malkus splits, for $\kdg\geq 2$, the inf-sup constant is bounded
and converges relatively rapidly in our tests, but for $\kdg=1$, the inf-sup constant 
goes to zero with a rate equal to one over the mesh size parameter~\cite{ref:QinThesis}.
The column ``iterations'' lists the number of iterations of the power method
required to achieve a change in $\lambda$ less than $10^{-8}$.
The column ``restarts'' lists the number of times the projection $\zz_\ell$ was
subtracted from $\uu^\pmi$.

In \cref{tabl:cubqwode}, we give results for Lagrange elements of degree $\kdg=2,3,4$
in two dimensions on Type I meshes, that is, meshes consisting of
45-degree right triangles (\cref{fig:twotypes}).
The results match the theory indicated in \cref{tabl:divfreeOK}: the inf-sup constant
degenerates for $\kdg=2,3$, but does not for $\kdg = 4$.
There is renewed interest in the low order cases due to the emergence of the grad-div penalized
Taylor--Hood method~\cite{case2011connection,linke2011convergence,ref:gradivtaylorhoodtosv}.
The approximations from these methods tend to divergence free functions as the penalty is increased,
and the low-order Taylor--Hood methods are widely used.
It has been known that divergence-free quadratics on Type I meshes have reduced approximation order
(this is equivalent to reduced approximation for $C^1$ piecewise cubic approximation for scalar functions~\cite{ref:seeonecubics}).
But so far we know of no estimates for controlled degeneration of the inf-sup constants
on general meshes.

\subsection{3D tests}
\label{sec:threond}

\begin{table}
\centering
\caption{Computation of inf-sup constants for various polynomial degrees $\kdg$
in three dimensions on the Freudenthal meshes depicted in \cref{fig:freudenthal-subdivision}.
The shift parameter in \cref{eqn:shtrfmsv} was set to $\shift=0.6$.
The mesh size $N$ refers to the $N \times N \times N$ mesh of cubes before subdivision.
Iterations were continued until the change in eigenvalue was less than $10^{-7}$
for $k=4$ and $k=5$, but we had to continue iterations 
until the change in eigenvalue was less than $10^{-8}$ for $k=3$.
\label{tabl:threedee}}
\begin{tabular}[t]{c|cccc}
\toprule
degree $\kdg$ & $N$ & iterations &inf-sup $\lambda_1$ &  restarts \\
\midrule
3& 2  & 1026   & $5.75 \times 10^{-4}$ & 4   \\
3& 3  & 3400   & $4.26 \times 10^{-4}$ & 14   \\
3& 4  & 4445   & $2.93 \times 10^{-4}$ & 29   \\
3& 5  & 6044   & $2.11 \times 10^{-4}$ & 2048   \\
3& 6  & 7732   & $1.59 \times 10^{-4}$ & 4863   \\
\midrule
4& 2  &  802   & $4.28 \times 10^{-3}$ & 7   \\
4& 3  & 2067   & $4.13 \times 10^{-3}$ & 16   \\
4& 4  & 1961   & $4.22 \times 10^{-3}$ & 16   \\
4& 5  & 1325   & $4.27 \times 10^{-3}$ & 12   \\
4& 6  & 1380   & $4.26 \times 10^{-3}$ & 12   \\
\midrule
5& 2  & 1004   & $6.46 \times 10^{-3}$ & 12   \\
5& 3  & 1003   & $6.59 \times 10^{-3}$ & 13   \\
5& 4  &  877   & $6.34 \times 10^{-3}$ & 12   \\
5& 5  &  938   & $6.26 \times 10^{-3}$ & 13   \\
5& 6  &  978   & $6.23 \times 10^{-3}$ & 13   \\
\bottomrule
\end{tabular}
\end{table}

Mesh restrictions in three dimensions for Scott--Vogelius elements are not
fully understood. Our computations for Freudenthal meshes are summarized in \cref{tabl:threedee}.
They reveal a surprising fact: on this mesh family,
the inf-sup constant is bounded for degrees $\kdg\geq 4$, whereas 
Zhang~\cite{zhang2011divergence} was able to prove inf-sup stability only for $\kdg\geq 6$.
Neilan~\cite[Proposition 6.5]{neilan2015discrete} extends this to more general meshes that 
satisfy a special condition, which he states as a conjecture, but still for $\kdg\geq 6$.
The results in \cref{tabl:threedee} lead us to make \cref{con:conjectureone}.

\subsection{Size of \texorpdfstring{{$\sdiv V^k_h$}}{{div V_h^k}}}

In general, $\sdiv V^k_h\subset DG^{k-1}_h$, the latter space being all
discontinuous piecewise polynomials of degree $k-1$.
In two dimensions, these spaces can be very close in size, 
differing because of singular vertices and the mean zero constraint
due to the homogeneous boundary conditions on $V^k_h$.
But in three dimensions, much less is known.
The difference between $\sdiv V^k_h$ and $DG^{k-1}_h$ is known to be quite large
for Freudenthal meshes.
The constraints on the latter space required to be satisfied to be in the former
space are listed in~\cite[page 691]{zhang2011divergence}.
For a single cube subdivided by 6 tetrahedra, the dimension of the quotient
space is 67 for $k=6$.
We should note that there is a typographical error in equation (3.67) in~\cite{zhang2011divergence},
which should read instead~\cite{zhangemail}
$$
\dim P_h=n^3(k+2)(k+1)k-3kn(n^2+n+2)+5.
$$
Note that $\dim DG^{k-1}_h=n^3(k+2)(k+1)k$, so the number of constraints
(the dimension of the quotient space) is $3kn(n^2+n+2)-5$.
The algorithm proposed by Rognes~\cite{ref:RognesAutoInfSupConst} can be used to compute the
number of constraints on general meshes.

\subsection{More general elements}

The algorithm described and tested here can be used more generally
for computing inf-sup constants.
The only restriction is that $\Pi_h=\sdiv V_h$.

\section{Space decompositions for the kernel of the divergence}

Our second conjecture relates to multigrid solvers for \cref{eq:linear-elasticity}. In the nearly
incompressible regime, the equation becomes nearly singular, and standard multigrid methods break down.
A key breakthrough for such problems was made by Sch\"oberl~\cite{schoberl1999}, who devised conditions on the relaxation and prolongation operators that guarantee that a multigrid method is parameter-robust in the nearly-singular regime.

\subsection{Background}

The key condition for the relaxation is best stated in terms of space decompositions and subspace corrections~\cite{xu1992}. The multigrid relaxation method we employ will be induced by a space decomposition
\begin{equation} \label{eq:space_decomposition}
V_h = \sum_{i=1}^J V_i
\end{equation}
where an equation for an approximation to the error is solved over each subspace $V_i$ and the updates combined additively or multiplicatively. Each $V_i$ is assumed small enough so that direct solvers can be afforded. For example, if $V_h = \operatorname{span}\{\phi_1, \phi_2, \dots, \phi_N\}$, and each $V_i$ is chosen to be $V_i = \operatorname{span}\{\phi_i\}$, then combining the updates additively would yield the Jacobi relaxation, while combining multiplicatively would yield the Gauss--Seidel relaxation. In a domain decomposition approach, each $V_i$ could be taken as the functions supported on one subdomain of a given parallel decomposition of the mesh, combined with a suitable global coarse space. Another important example is to define each $V_i$ as the functions with support on the patch of cells surrounding each vertex, the so-called \emph{vertex-star} space decomposition: this arises in $k$-robust preconditioners for symmetric and coercive problems~\cite{pavarino1993,schoberl2008}, as the Arnold--Falk--Winther (AFW) relaxation for $H(\div)$ and $H(\curl)$~\cite{arnold2000}, and in Reynolds-robust preconditioners for the Navier--Stokes equations~\cite{benzi2006,farrell2018b}.

Let us consider an abstract nearly singular problem, following Lee et al.~\cite{lee2007}. For $\varepsilon > 0$, consider the finite-dimensional linear variational problem: find $u \in V_h$ such that
\begin{equation} \label{eq:generic_nearly_singular}
a_0(u, v) + \varepsilon^{-1} a_1(u, v) = (f, v) \text{ for all } v \in V_h,
\end{equation}
where $a_0$ is symmetric and coercive, and $a_1$ is symmetric but only positive semi-definite. In the context of \cref{eq:linear-elasticity}, $a_1(u, v) = (\nabla \cdot u, \nabla \cdot v)_{L^2}$. Define the kernel
\begin{equation}
\mathcal{N} = \{u \in V_h: a_1(u, v) = 0 \text{ for all } v \in V_h\}.
\end{equation}
In our context, these are divergence-free functions in the finite element space. A key condition that must be satisfied by the space decomposition for parameter-robustness in $\varepsilon$ is~\cite[Assumption (A1)]{lee2007}:
\begin{equation} \label{eq:assumption_A1}
\mathcal{N} = \sum_{i=1}^J \left( \mathcal{N} \cap V_i \right).
\end{equation}
In other words, when challenged with a divergence-free function, it must be possible to decompose this as the sum of functions with the $i^{\mathrm{th}}$ summand drawn from the divergence-free functions in $V_i$\footnote{This decomposition must also be stable, but we shall not elaborate here.}. The space decomposition that induces Jacobi or Gauss--Seidel relaxation does not generally satisfy \cref{eq:assumption_A1}, because the intersection $\mathcal{N} \cap V_i = \{0\}$ for typical finite element methods where the standard bases used are not divergence-free.

\subsection{Space decompositions from de Rham complexes}
One way to devise space decompositions that satisfy \cref{eq:assumption_A1} is by inspecting discrete subcomplexes of a suitable underlying Hilbert complex. For concreteness, consider $\Omega \subset \mathbb{R}^2$, with $\Omega$ simply connected. The Stokes complex~\cite{john2017} is given by
\begin{align}
  \label{eq:stokes_complex_2d}
  \mathbb{R} \xrightarrow{\operatorname{id}} H^2(\Omega)
  \xrightarrow{\curl} [H^1(\Omega)]^2 \xrightarrow{\div} L^2(\Omega) \xrightarrow{\operatorname{null}} 0.
\end{align}
This complex is discretized with a discrete subcomplex
\begin{equation} \label{eq:discrete_stokes_complex_2d}
 \mathbb{R} \xrightarrow{\operatorname{id}} \Sigma_h \xrightarrow{\curl} V_h \xrightarrow{\div} Q_h \xrightarrow{\operatorname{null}} 0,
\end{equation}
where $\Sigma_h \subset H^2(\Omega)$, $V_h \subset [H^1(\Omega)]^2$, and $Q_h \subset L^2(\Omega)$. These complexes have the property that the kernel of an operator $\curl, \div$, or $\operatorname{null}$ is a subspace of the range of the preceding operator, e.g.~that $\operatorname{ker}(\div) \subset \operatorname{range}(\curl)$~\cite{arnold2018}. The complex is called \emph{exact} if the kernel of an operator is precisely the range of the preceding operator. The Stokes complex is exact if the domain $\Omega$ is simply connected, and this property is inherited by \cref{eq:discrete_stokes_complex_2d} if constructed appropriately (with bounded cochain projections).

In our context, these complexes are useful because they offer a crisp characterization of $\mathcal{N} = \operatorname{ker}(\div)$. Let $\Sigma_h = \operatorname{span}\{\psi_1, \dots, \psi_M\}$ and $V_h = \operatorname{span}\{\phi_1, \dots, \phi_N\}$. If $u_h \in \mathcal{N}$, then $u_h = \curl \Psi_h$ for $\Psi_h \in \Sigma_h$ by exactness of \cref{eq:discrete_stokes_complex_2d}. Expanding $\Psi_h$ in terms of its basis functions
\begin{equation}
\Psi_h = \sum_{i=1}^M c_i \psi_i
\end{equation}
yields an expression for $u_h$
\begin{equation}
u_h = \sum_{i=1}^M c_i \curl \psi_i.
\end{equation}
Suppose the space decomposition \cref{eq:space_decomposition} is chosen as $J=N + M$ with
\begin{equation} \label{eq:impractical_choice}
V_i = \begin{cases}
\operatorname{span}(\curl \psi_i) & i = 1, \dots, M \\
\operatorname{span}(\phi_{i-M})       & i = M+1, \dots, M+N.
\end{cases}
\end{equation}
This decomposition would satisfy \cref{eq:assumption_A1}, with $V_i \cap \mathcal{N} = V_i$ for $i \le M$ or $V_i \cap \mathcal{N} = \{0\}$ for $i > M$\footnote{This space decomposition relates to that proposed by Hiptmair~\cite{hiptmair1998}. As written it is not robust in $k$, but (experimentally) can be made so by applying suitable block Jacobi methods in $V_h$ and $\Sigma_h$.}.

Solving the problem \cref{eq:generic_nearly_singular} over $\operatorname{span}\{\curl \psi_i\}$ reduces to solving a problem in a locally-supported subspace of $\Sigma_h$, as described in Hiptmair~\cite{hiptmair1998} and Hiptmair--Xu~\cite{hiptmair2007} for the $L^2$ de Rham complex. However, for the Stokes complex in three dimensions, explicit constructions of $\Sigma_h$ are poorly understood. The alternative approach (which we refer to as Pavarino--Arnold--Falk--Winther, PAFW) is to construct the space decomposition using only knowledge of the supports of the basis functions $\psi_i$. For example, if $\Sigma_h$ exists with the property that each basis function $\psi_i$ is a polynomial of degree $k+1$ supported on a certain region of the domain $\operatorname{supp}(\psi_i)$, then we may choose $\{V_j\}$ so that each $V_j$ captures all polynomials of degree $k$ on one $\operatorname{supp}(\psi_i)$. In this way, $\forall i =1, \dots, M\, \exists\, j \text{ s.t. } \curl \psi_i \in V_j$, and \cref{eq:assumption_A1} will also be satisfied, without needing to explicitly construct $\Sigma_h$ or $\{\psi_i\}_{i=1}^M$. Only knowledge (or conjecture) of the supports is required.

To illustrate the PAFW construction, consider the case of the Scott--Vogelius element pair with $k \ge 4$ and $d = 2$. This element constructs $V_h$ with $[\mathrm{CG}_k]^d$ and $Q_h$ with $\mathrm{DG}_{k-1}$ on a mesh $\mathcal{T}_h$. The corresponding $C^1$-conforming finite element for $\Sigma_h$ is the Morgan--Scott element of degree $k+1$~\cite{lrsBIBaf}, which is not implemented in general purpose finite element software, preventing the application of \cref{eq:impractical_choice}. The Morgan--Scott element employs (among others) degrees of freedom at the vertices of the mesh, and the basis function associated with such a degree of freedom will have support over the patch of cells sharing that vertex. This suggests the following vertex-star space decomposition, which we now describe. Let $\vt_1, \dots, \vt_J$ be the vertices of $\mathcal{T}_h$. The subspaces are given by
\begin{equation} \label{eq:vertex_star}
V_j = \{v_h \in V_h : \operatorname{supp}(v_h) \subset \operatorname{star}(\vt_j)\},
\end{equation}
where the $\operatorname{star}$ operation of a simplex $p$ returns the union of all simplices containing $p$ as a subsimplex~\cite[\S 2]{munkres1984}. In other words,
when applied to a vertex $\vt_j$, the star is the union of the cells and edges (and faces in $d=3$) containing $\vt_j$, as well as $\vt_j$ itself. 
In~\cite{lrsBIBaf}, it is proven that the Morgan--Scott basis functions are each supported in
$\operatorname{star}(\vt_j)$ for some $j$, provided $k\geq 4$.
Thus \cref{eq:vertex_star} also satisfies the kernel decomposition \cref{eq:assumption_A1}. Numerical experiments indicate that \cref{eq:vertex_star} does indeed yield $\varepsilon$-robust convergence for $k \ge 4$, as expected from the Morgan--Scott theory, while it does not for $k < 4$~\cite[\S 4.2]{farrell2019c}.

There are other examples where the vertex-star space decomposition does not yield $\varepsilon$-robust convergence. For $k=1$ on Malkus splits, the basis with smallest possible support has support
that is larger than $\operatorname{star}(\vt_j)$, although it is in $\operatorname{star}(\operatorname{closure}(\operatorname{star}(\vt_j)))$~\cite{lrsBIBim}, where the $\operatorname{star}$ of a set of simplices is the union of their stars, and closure is defined in~\cite{farrell2019c}.

\subsection{Second conjecture}

The preceding discussion indicates that one may experimentally
investigate whether $\Sigma_h$ permits a basis with support
captured by the stars of the vertices, by applying the multigrid
algorithm and observing whether the convergence is $\varepsilon$-robust
or not.

The continuous Stokes complexes in 3D~\cite{farrell2020reynolds} is given by
\begin{equation} \label{eqn:stokescomplex3dcts}
  \mathbb{R} \xrightarrow{\operatorname{id}} H^2(\Omega)
 \xrightarrow{\grad} H^1(\curl, \Omega)
  \xrightarrow{\curl} [H^1(\Omega)]^3
  \xrightarrow{\div} L^2(\Omega)
  \xrightarrow{\operatorname{null}} 0.
\end{equation}
Here
\begin{equation}
H^1(\curl, \Omega) \coloneqq \{v \in [H^1(\Omega)]^3: \curl v \in [H^1(\Omega)]^3\}.
\end{equation}
This complex is discretized with a discrete subcomplex
\begin{equation} \label{eqn:stoplextred}
  \mathbb{R} \xrightarrow{\operatorname{id}} S_h
 \xrightarrow{\grad} \Sigma_{h}
  \xrightarrow{\curl} V_{h}
  \xrightarrow{\div} \Pi_{h}
  \xrightarrow{\operatorname{null}} 0.
\end{equation}
The space $S_h$ consists of $C^1$ piecewise polynomials of degree $k+2$.
The space $\Sigma_{h}$ consists of continuous vector-valued piecewise polynomials of degree $k+1$
which have a continuous curl.
Note that $(S_h)^3\subset \Sigma_{h}$.

The potential space $\Sigma_h$ is known on special meshes (Alfeld or Worsey--Farin splits) in three dimensions~\cite{fu2020exact,guzman2020exact,boffi2022convergence}.
We are not aware of any results that characterize $\Sigma_h$
on Freudenthal meshes.

Based on the numerical experiments we will report, we make the
following conjecture:
\begin{conjecture} \label{con:conjecturetwo}
Let $V_h \subset [H^1(\Omega)]^3$ be constructed with
continuous Lagrange elements of degree $k$,
and let $\Sigma_h \subset H^1(\curl, \Omega)$
be the space preceding $V_h$ in
a subcomplex of the three-dimensional Stokes complex.
\begin{conjecturelist}[label={\normalfont (\alph*)}, ref={\theconjecture(\alph*)}]
\item \label{con:conjecturetwo-a} For $\kdg = 4$, there does not exist a local basis for $\Sigma_h$, supported 
on the stars of vertices, on the Freudenthal meshes depicted in 
\cref{fig:freudenthal-subdivision}.
\item \label{con:conjecturetwo-b} For $\kdg \ge 5$, there exists a local basis for $\Sigma_h$, supported on
the stars of vertices, on the Freudenthal meshes depicted in 
\cref{fig:freudenthal-subdivision}.
\end{conjecturelist}
\end{conjecture}

Our evidence for \cref{con:conjecturetwo-b} is not as strong as for \cref{con:conjecturetwo-a}.
For example, the results of Sch\"oberl~\cite{schoberl1999} are not
necessary and sufficient, so robustness could occur for some other reason.
But lack of robustness implies lack of a suitable local basis.

\subsection{Numerical experiments}

\begin{figure}[tbhp]
\footnotesize
\centering
\begin{tikzpicture}[%
 every node/.style={draw=black, thick, anchor=west},
grow via three points={one child at (-0.0,-0.7) and
	two children at (0.0,-0.7) and (0.0,-1.4)},
edge from parent path={(\tikzparentnode.210) |- (\tikzchildnode.west)}]
\node {Krylov solver: CG}
child {node {Two-grid V-cycle}
   child {node {Relaxation: damped additive vertex-star iteration}
   }
   child {node {Coarse grid: Cholesky factorization}
   }
};
\end{tikzpicture}
\caption{Solver diagram for \cref{eq:nearly_singular}.}
\label{fig:solver}
\end{figure}
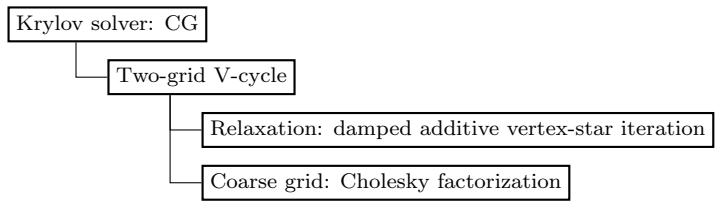

We provide numerical evidence for \cref{con:conjecturetwo}.
We solve the problem:
find $u \in V_h \subset [{H}^1_0(\Omega)]^3$ such that:
\begin{equation} \label{eq:nearly_singular}
  (\nabla u , \nabla v)_{L^2(\Omega)} + \gamma (\sdiv u, \sdiv v)_{L^2(\Omega)}
  = (1, v)_{L^2(\Omega)} \quad \forall v \in V_h.
\end{equation}
The problem becomes nearly singular as $\gamma \to \infty$.
Here $V_h$ is constructed with continuous Lagrange elements of degree $\kdg$. We employ the solver denoted in \cref{fig:solver}.
The relaxation on each level is one application of damped Richardson
iteration preconditioned by the additive Schwarz method with vertex-star subspaces chosen as \cref{eq:vertex_star}.
The damping factor is set to $1/M$ where $M$ is the maximum number 
of patches that a given vertex is contained in. Here $M = 3$ for $d = 2$ and $M = 4$ for $d = 3$.
The coarse grid solve was computed with CHOLMOD~\cite{chen2008} via PETSc~\cite{balay2022}.
The code was run in serial.

\begin{table}[htbp]
  \centering
\caption{Iteration counts for \cref{eq:nearly_singular} on a uniform
refinement of a Type I $4 \times 4$ mesh. The iteration counts are $\gamma$-robust
for $k \ge 4$, but not for $k < 4$.\label{tab:two-dim-richardson}}
\begin{tabular}{l|cccccc}
  \toprule
 $\kdg\backslash \gamma$ & $10^0$ & $10^1$ & $10^2$ & $10^3$ & $10^4$ & $10^5$ \\
 \midrule
 2                    & 11     & 17     & 35     & 54     & 87     & 100    \\
 3                    & 10     & 13     & 23     & 49     & 62     & 83     \\
 4                    & 9      & 10     & 13     & 13     & 13     & 12     \\
 5                    & 9      & 10     & 11     & 12     & 11     & 10     \\
  \bottomrule
  \end{tabular}
\end{table}

\subsection{Two dimensions}

The convergence results for the solver depicted in \cref{fig:solver} are shown in \cref{tab:two-dim-richardson}, as a function of $\kdg$ and
$\gamma$. The Krylov method was terminated when the unpreconditioned residual was reduced by 8 orders of magnitude.
The coarse mesh was a Type I mesh~\cite{lai2007spline} of size $4\times 4$.
The main observation is that $\kdg\geq 4$ is required
for $\gamma$-robust performance.
This matches well with the existing theory described above: the key property that changes between these two cases is the
existence of the Morgan--Scott element for $\Sigma_h$ with a basis supported on vertex-stars for $k \ge 4$.

\subsection{Three dimensions}

In three dimensions we again choose $V_h = [\mathrm{CG}_k]^d$.
We employ an analog of a Type I coarse mesh of size $4 \times 4 \times 4$, with one refinement.
We show iteration counts for the solution of
the problem of \cref{eq:nearly_singular} in \cref{tab:three-dim-richardson}.
As in two dimensions, the solver was terminated when the
unpreconditioned residual was reduced by 8 orders of magnitude. These
results lead us to \cref{con:conjecturetwo}.
\begin{table}[htbp]
\centering
\caption{Iteration counts for \cref{eq:nearly_singular} on 2 uniform
refinements of a Type I $4 \times 4 \times 4$ mesh. The iteration counts are $\gamma$-robust
for $k \ge 5$, but not for $k < 5$.\label{tab:three-dim-richardson}}
\begin{tabular}{l|cccccc}
  \toprule
    $\kdg\backslash \gamma$ & $10^0$ & $10^1$ & $10^2$ & $10^3$ & $10^4$ & $10^5$ \\
    \midrule
    2                    & 14     & 23     & 52     & 118    & 300    & 849    \\
    3                    & 11     & 16     & 29     & 66     & 173    & 458    \\
    4                    & 11     & 13     & 19     & 26     & 54     & 110    \\
    5                    & 11     & 12     & 16     & 19     & 20     & 19     \\
    6                    & 10     & 11     & 14     & 15     & 16     & 15     \\
    7                    & 10     & 11     & 13     & 14     & 14     & 13     \\
  \bottomrule
  \end{tabular}
\end{table}

Note that our \cref{con:conjectureone} is that the global inf-sup
condition holds on this mesh family for $\kdg\geq 4$.
Thus there is a gap for $\kdg=4$: we conjecture that the global inf-sup 
condition holds, but that there is not a basis for the potential space
supported on vertex-stars.

\section{Conclusions}

Our computational experiments have suggested two conjectures regarding the
three-dimensional Stokes complex on Freudenthal meshes.
We conjecture that the inf-sup condition holds for a velocity space $V_h$ consisting
of continuous piecewise polynomials of degree $\kdg$ and pressure space $\sdiv V_h$
for $\kdg \ge 4$,
and that there is a local basis for the associated potential space $\Sigma_h$
for $\kdg \ge 5$.

\section{Acknowledgements}

The authors would like to thank  Michael Neilan, Charles Parker, Florian Wechsung,
Umberto Zerbinati, and Shangyou Zhang for useful discussions.

\appendix

\section{Freudenthal triangulation of the unit cube} \label{sec:freudenthal}

For concreteness we describe the Freudenthal triangulation of the unit cube $[0, 1]^3$. Define vertices with coordinates
\begin{align*}
\textbf{v}_1 = (0, 0, 0),\quad & \textbf{v}_2 = (1, 0, 0), \quad
\textbf{v}_3 = (1, 1, 0),\quad  \textbf{v}_4 = (1, 1, 1), \\
\textbf{v}_5 = (0, 1, 0),\quad & \textbf{v}_6 = (1, 0, 1), \quad
\textbf{v}_7 = (0, 1, 1),\quad  \textbf{v}_8 = (0, 0, 1).
\end{align*}
Then the six cells (tetrahedra) $\textbf{c}_i$ have vertices
\begin{align*}
\textbf{c}_1 = [\textbf{v}_1, \textbf{v}_2, \textbf{v}_3, \textbf{v}_4], \quad & 
\textbf{c}_2 = [\textbf{v}_1, \textbf{v}_3, \textbf{v}_4, \textbf{v}_5],\quad
\textbf{c}_3 = [\textbf{v}_1, \textbf{v}_2, \textbf{v}_4, \textbf{v}_6],\\ 
\textbf{c}_4 = [\textbf{v}_1, \textbf{v}_4, \textbf{v}_5, \textbf{v}_7],\quad & 
\textbf{c}_5 = [\textbf{v}_1, \textbf{v}_4, \textbf{v}_6, \textbf{v}_8],\quad 
\textbf{c}_6 = [\textbf{v}_1, \textbf{v}_4, \textbf{v}_7, \textbf{v}_8].
\end{align*}
It is well known~\cite[Figure 2]{zhang2011divergence} that a standard multigrid
subdivision of a Freudenthal triangulation of the unit cube $[0, 1]^3$ yields
eight subcubes, each subdivided by a Freudenthal triangulation.
In particular, the subdivision of the Freudenthal triangulation of the unit cube 
$[0, 1]^3$ consists of cutting by the three planes $x=1/2$, $y=1/2$, and $z=1/2$.
See~\cite[section 4]{zhang2011divergence} for more information regarding 
multigrid for these meshes.

\ifarxiv

\input{paper.bbl}
\else
\bibliographystyle{siamplain}
\bibliography{paper}
\fi
\end{document}

%% file: freudenthal-subdivision.tex
\begin{tikzpicture}[3d view={120}{10}, scale=2]
  \coordinate [label=below:$\textbf{v}_1$] (v1) at (0, 0, 0) {};
  \coordinate [label=below:$\textbf{v}_2$] (v2) at (1, 0, 0) {};
  \coordinate [label=below:$\textbf{v}_3$] (v3) at (1, 1, 0) {};
  \coordinate [label=above:$\textbf{v}_4$] (v4) at (1, 1, 1) {};
  \coordinate [label=below:$\textbf{v}_5$] (v5) at (0, 1, 0) {};
  \coordinate [label=above:$\textbf{v}_6$] (v6) at (1, 0, 1) {};
  \coordinate [label=above:$\textbf{v}_7$] (v7) at (0, 1, 1) {};
  \coordinate [label=above:$\textbf{v}_8$] (v8) at (0, 0, 1) {};

  \begin{scope}[xshift=-1.25cm, yshift=0.125cm]
    \draw[-stealth] (0, 0, 0) -- (0.5, 0, 0) node[pos=1.2] {x};
    \draw[-stealth] (0, 0, 0) -- (0, 0.5, 0) node[pos=1.1] {y};
    \draw[-stealth] (0, 0, 0) -- (0, 0, 0.5) node[pos=1.1] {z};
  \end{scope}
  \draw (v1) -- (v2) (v1) -- (v3) (v1) -- (v4);
  \draw (v2) -- (v3) (v2) -- (v4);
  \draw (v3) -- (v4);

  \draw (v1) -- (v3) (v1) -- (v4) (v1) -- (v5);
  \draw (v3) -- (v4) (v3) -- (v5);
  \draw (v4) -- (v5);

  \draw (v1) -- (v2) (v1) -- (v4) (v1) -- (v6);
  \draw (v2) -- (v4) (v2) -- (v6);
  \draw (v4) -- (v6);

  \draw (v1) -- (v4) (v1) -- (v5) (v1) -- (v7);
  \draw (v4) -- (v5) (v4) -- (v7);
  \draw (v5) -- (v7);

  \draw (v1) -- (v4) (v1) -- (v6) (v1) -- (v8);
  \draw (v4) -- (v6) (v4) -- (v8);
  \draw (v6) -- (v8);

  \draw (v1) -- (v4) (v1) -- (v7) (v1) -- (v8);
  \draw (v4) -- (v7) (v4) -- (v8);
  \draw (v7) -- (v8);
\end{tikzpicture}

%% file: twotypes.tex
\begin{tikzpicture}
  \draw (0, 0) -- (2, 0) -- (2, 2) -- (0, 2) -- cycle;
  \draw (0, 0) -- (2, 2) (0, 1) -- (2, 1) (1, 2) -- (1, 0);
  \draw (0, 1) -- (1, 2) (1, 0) -- (2, 1);

  \node at (1, -0.5) {Type I};
  \begin{scope}[xshift=3cm]
  \draw (0, 0) -- (2, 0) -- (2, 2) -- (0, 2) -- cycle;
  \draw (0, 1) -- (2, 1) (1, 0) -- (1, 2);
  \draw (0, 0) -- (2, 2) (0, 2) -- (2, 0);
  \draw (0, 1) -- (1, 2) (1, 0) -- (2, 1);
  \draw (0, 1) -- (1, 0) (1, 2) -- (2, 1);
  \node at (1, -0.5) {Type II};
  \end{scope}
\end{tikzpicture}